\documentclass{article} 
\usepackage{amssymb,amsthm,amsmath,amsfonts}
\usepackage{graphicx}
\usepackage{cite}
  \newtheorem{theorem}{Theorem}

               \newtheorem{lemma}{Lemma}
               
               \newtheorem{remark}{Remark}
               
               \def\pf{\par\noindent {\em Proof.}~\par\noindent}
               \def\qed{~\hfill{$\square$}\pagebreak[1]\par\medskip\par}

\newcommand{\f}{{\bf{f}}}
\newcommand{\Li}{{\mbox{Lip}}}
\newcommand{\bl}{{\bf l}}
\newcommand{\bu}{{\bf u}}

\newcommand{\bg}{{\bf g}}
\newcommand{\bj}{{\bf j}}

\newcommand{\cQ}{{\cal Q}}

\newcommand{\R}{{\mathbb R}}

\newcommand{\ux}{\underline{x}}
\newcommand{\pD}{{^\psi\!\upa}}
\newcommand{\hD}{{^\varphi\!\upa}}

\newcommand{\uy}{\underline{y}}
\newcommand{\uz}{\underline{z}}

\newcommand{\upa}{\underline{\partial}}

\newcommand{\Om}{\Omega}
\newcommand{\cH}{{\cal H}}
\newcommand{\cW}{{\cal W}}

\newcommand{\cC}{{\cal C}}

\newcommand{\cT}{{\cal T}}

\newcommand{\pa}{\partial}

\begin{document}
\title{On a Riemann-Hilbert boundary value problem for $(\varphi,\psi)$-harmonic functions in $\R^m$}
\small{
\author
{Jos\'e Luis Serrano Ricardo$^{(1)}$; Ricardo Abreu Blaya$^{(1)}$;\\Juan Bory Reyes$^{(2)}$; Jorge S\'anchez Ortiz$^{(1)}$}
\vskip 1truecm
\date{\small $^{(1)}$ Facultad de Matem\'aticas, Universidad Aut\'onoma de Guerrero, M\'exico.\\ Emails: jose.luis.serrano.999@gmail.com, rabreublaya@yahoo.es, jsanchezmate@gmail.com\\$^{(2)}$ {SEPI-ESIME-ZAC-Instituto Polit\'ecnico Nacional, Ciudad M\'exico, M\'exico}\\Email: juanboryreyes@yahoo.com}

\maketitle
\begin{abstract}
The purpose of this paper is to solve a kind of Riemann-Hilbert boundary value problem for $(\varphi,\psi)$-harmonic functions, which are linked with the use of two orthogonal basis of the Euclidean space $\R^m$. We approach this problem using the language of Clifford analysis for obtaining the explicit expression of the solution of the problem in a Jordan domain $\Omega\subset\R^m$ with fractal boundary. One of the remarkable feature in this study is that the boundary data involves higher order Lipschitz class of functions.
\end{abstract}

\vspace{0.3cm}

\small{
}
\section{Introduction}
We consider the $2^m$-dimensional real Clifford algebra $\R_{0,m}$ generated by $e_1,e_2,\ldots,e_m$ according to the multiplication rules 
$$e_ie_j+e_je_i=-2\delta_{i,j},$$ 
where $\delta_{i,j}$ is the Kronecker's symbol.

The elements of the algebra $\R_{0,m}$ have a unique representation of the form
$$a=\sum_{A\subseteq \{1,2,\ldots,m\}}a_Ae_A,$$
where $a_A\in\R$ and where we identify $e_A$ with $e_{h_1}e_{h_2}\cdots e_{h_k}$ for $A=\{h_1, h_2, \ldots,h_k\}$ $(1\leq h_1<h_2<\cdots<h_k\leq m)$ and $e_\emptyset=e_0=1$.

We can embed $\R^m$ into $\R_{0,m}$ by identifying $\ux = (x_1, . . . , x_m)$ with $\ux = \sum_{i=1}^m x_ie_i$. All properties necessary of this algebras one can find in \cite{GM} for instance.

Clifford analysis is nowadays described as a theory of $\R_{0,m}$-valued functions of $m$ real variables based on the null-solutions of the Dirac operator (a generalization of the complex Cauchy-Riemann operator) offering a refinement of classical harmonic analysis. For more information regarding Clifford analysis, see, e.g., \cite{BDS, GHS, G16}. 

Suppose $\Omega\subset\R^{m}$ is a Jordan domain, with boundary a compact topological surface $\Gamma$, which decompose $\R^m$ into the interior and exterior (containing the infinity point) domains denoted by $\Omega_+$ and $\Omega_-$ respectively. Under further notice we assume $\Gamma$ to be sufficiently smooth.

We will be interested in functions $u:\Om\rightarrow$ $\R_{0,m}$, which might be written as $u(\ux)=\sum_{A}u_A(\ux)e_A$ with $u_A$ $\R$-valued. Properties, such continuity, differentiability, integrability, and so on, are ascribed coordinate-wise. In particular, we define in this way the following right module of $\R_{0,m}$-valued functions:
\begin{itemize}
\item $C^k(\Omega,\R_{0,m}), k\in\mathbb{N}\cup\{0\}$ the right module of all $\R_{0,m}$-valued functions, $k$-times continuously differentiable in $\Omega$. For $u\in C^k(\Omega,\R_{0,m})$ we will write 
$$\partial^{\bj}=\frac{\partial^{|\bj|}}{\partial x_1^{j_1}\partial x_2^{j_2}\dots \partial x_m^{j_m}},$$
where $\bj= (j_1,j_2,\dots , j_m)\in\left(\mathbb N\cup\{0\}\right)^{m}$ is a $m$-dimensional multi-indices and $|\bj|=j_1+\cdots+j_m$.
\item $C^{k,\alpha}(\Omega,\R_{0,m})$, $\alpha \in (0, 1]$ the right module of all $\R_{0,m}$-valued functions, $k$-times $\alpha$-H\"older continuously differentiable in $\Omega$.
\item $L_p(\Omega,\R_{0,m})$, ($1\leq p<\infty$) the right module of all equivalence classes of $p$-Lebesgue measurable $\R_{0,m}$-functions over $\Omega$.
\end{itemize}
Let an ordered set $\psi:=\{\psi^1,\ldots,\psi^m\}$, with $\psi^i\in\R^m\subset\R_{0,m}$. On the set $C^1(\Omega,\R_{0,m})$ we define the generalized Dirac operator by: 
\begin{equation}\label{DiracPsi}
\pD:=\psi^1{\frac{\partial}{\partial x_1}}+\psi^2{\frac{\partial}{\partial x_2}}+\cdots\psi^m{\frac{\partial}{\partial x_m}}.
\end{equation}
For the particular case of the standard $\R_{0,m}$-basic vector set
$\psi_{st}:=\{e_1, e_2\ldots,e_m\},$
operator $\pD$ becomes the Dirac operator.

Let $\Delta_{m}$ be the $(m)$-dimensional Laplace operator. It is easy to prove that the equality
\begin{equation}\label{Fac Lap}
\pD\pD=-\Delta_{m} 
\end{equation}
in $\R^m$ hold, if and only if 
\[\psi^i\psi^j+\psi^j\psi^i=-2\delta_{ij}\,\,(i,j=1,2,\dots m).\]
Note that last equality yields
\begin{equation}\label{Estr Cond}
2\delta_{i,j}=\psi^i\cdot \bar{\psi^j}+ \psi^j\cdot\bar{\psi^i}=2\left\langle \psi^i, \psi^j\right\rangle_{\R^{m}},
\end{equation}
where $\left\langle , \right\rangle_{\R^{m}}$ denotes the scalar product, hence factorization (\ref{Fac Lap}) holds if and only if $\psi$ represents an orthonormal basis of $\R^{m}$. 

A set $\psi$ with the property (\ref{Estr Cond}) is called {\it structural set}. Notion of structural sets goes back to \cite{No, S}.

The $\R_{0,m}$-valued solutions of $\pD u=0$ are the so-called $\psi$-hyperholomorphic functions. Next we mention basic facts of this theory to be used in the paper, thus making our exposition self-contained. Deeper discussions can be found in \cite{A1, A2, BDGS} and the references given there.

The fundamental solution of the operator $\pD$ is given by
\[
K_\psi(\ux)=\frac{-\ux_\psi}{\sigma_m|\ux|^m},
\]
where
\[
\ux_\psi=\sum_{i=1}^m x_i\psi^i\,\,\,\mbox{if}\,\,\,\, \ux=\sum_{i=1}^m x_i e_i
\]
and $\sigma_m$ stands for the area of the unit sphere in $\R^m$.

This particularly important function, referred to as Cauchy kernel, comes from acting $\pD$ to the fundamental solution of the Laplacian $\Delta_{m}$ given by $\displaystyle\frac{|\ux|^{2-m}}{\sigma_m(2-m)}$, i.e.
\[
K_\psi(\ux)=\pD[\frac{|\ux|^{2-m}}{\sigma_m(2-m)}].
\]
The Cauchy kernel $K_\psi$ is $\psi$-hyperholomorphic in $\R^m\setminus\{0\}$ and plays a decisive role in our context. 
\begin{theorem}[Borel-Pompeiu formula]
Let $u\in C^1(\Omega\cup \Gamma,\R_{0,m})$. Then it holds that
\begin{equation}\label{BP0}
\int_{\Gamma}K_\psi(\uy-\ux)n_\psi(\uy)u(\uy)d\uy-\int_\Omega K_\psi(\uy-\ux)\pD u(\uy)d\uy=\begin{cases}u(\ux)&\,\mbox{if $\ux\in\Omega_+$}\\0&\,\mbox{if $\ux\in\Omega_-$,}\end{cases}
\end{equation}
\end{theorem}
\noindent
where $n_\psi(\uy)=\sum_{i=1}^m n_i(\uy)\psi^i$, being $n_i(\uy)$ the $i$-th component of the outward unit normal vector at $\uy\in\Gamma$.

From \eqref{BP0} one finds two important integral operators: the Cauchy transform
\[
\cC_\psi u(\ux):=\int_{\Gamma}K_\psi(\uy-\ux)n_\psi(\uy)u(\uy)d\uy,
\] 
which represents a $\psi$-hyperholomorphic function in $\R^m\setminus\Gamma$ and the Teodorescu operator
\[
\cT_\psi v(\ux)=-\int_\Omega K_\psi(\uy-\ux)v(\uy)d\uy,
\]
which runs as the right-handed inverse of $\pD$, i.e.
\begin{equation}\label{rel1}
\pD\cT_\psi v(\ux)=\begin{cases}v(\ux),&\mbox{if $\ux\in\Omega_+$}\\0&\mbox{if $\ux\in\Omega_-$}.\end{cases}
\end{equation}
\section{Preliminaries}
\subsection{Lipschitz classes and Whitney extension theorem}
The higher order Lipschitz class $\Li(k+\alpha,\Gamma)$ consists of the collections of real-valued continuous functions
\begin{equation}\label{LD}
\f:=\{f^{(\bj)},\,|\bj|\le k\}
\end{equation} 
defined on $\Gamma$ and satisfying the compatibility conditions
\begin{equation}\label{C}
|f^{(\bj)}(\ux)-\sum_{|\bj+\bl|\le k}\frac{f^{(\bj+\bl)}(\uy)}{\bl !}(\ux-\uy)^\bl|=\mathcal{O}(|\ux-\uy|^{k+\alpha-|\bj|}),\,\,\ux,\uy\in\Gamma,\,|\bj|\le k.
\end{equation}
In 1934 the American mathematician Hassler Whitney, one of the most prominent figures from the field in the 20th century \cite{K}, proved in \cite{Wh} that such a collection can be extended as a $C^{k,\alpha}$-smooth function on $\R^m$. For an excellent reference alone more classical lines we refer the reader to the well-known book of E. M. Stein \cite[Chapter VI, p. 176]{St}.

Because the $\R_{0,m}$-valued Lipschitz classes are component-wise defined, by abuse of notation we continue to write $\Li(k+\alpha,\Gamma)$ for the Cliffordian situation. We shall confine ourselves to discussing the case $k=1$, which will become clear shortly.

\begin{theorem}[Whitney]\label{Wh}
Let $\bu=\{u^{(\bj)},\,|\bj|\le 1\}$ be an $\R_{0,m}$-valued collection in $\Li(1+\alpha,\Gamma)$. Then, there exists a compact supported $\R_{0,m}$-valued function $\tilde{u}\in C^{1,\alpha}(\R^{m},\R_{0,m})$ satisfying
\begin{itemize}
\item[(i)] $\tilde{u}|_\Gamma = u^{(0)},\,\partial^{(\bj)}\tilde{u}|_\Gamma=u^{(\bj)},\,|\bj|=1$,
\item[(ii)] $\tilde{u}\in C^\infty(\R^{m} \setminus \Gamma)$,
\item[(iii)] $|\partial^{\bj}\tilde{u}(\ux)| \leqslant c \, \mbox{\em dist}(\ux,\Gamma)^{\alpha-1}$, for $|\bj|=2$ and $\ux\in\R^m\setminus\Gamma$.
\end{itemize}
\end{theorem}

\subsection{Harmonic versus $(\varphi,\psi)$-harmonic functions}
This paper is concerned with a second order partial differential equation arising in a natural way from the consideration of two different structural sets $\varphi$ and $\psi$, i.e. the generalized Laplace equation 
\begin{equation}\label{GL}
\hD\pD u=0.
\end{equation}
The solutions of \eqref{GL} will be referred as $(\varphi,\psi)$-harmonic functions. The space of all $(\varphi,\psi)$-harmonic functions in $\R^m$ will be denoted by $\cH_{\varphi,\psi}(\Omega,\R_{0,m})$.  It is worth noting that for $\varphi=\psi$, the class $\cH_{\varphi,\psi}(\Omega,\R_{0,m})$ coincides with the space $\cH(\Omega,\R_{0,m})$ of harmonic functions in $\Omega$, which justifies the name we choose for the functions in $\cH_{\varphi,\psi}(\Omega,\R_{0,m})$.

Further, it is likewise easy to obtain the fundamental solution of $\hD\pD$. In fact it is sufficient to apply the corresponding operator $\pD\hD$ to the fundamental solution of the bi-Laplacian $\Delta^{2}_{m}$, which is given by
\[
\displaystyle\frac{|\ux|^4}{\sigma_m (2-m)(4-m)}.
\]
This gives rise to the fundamental solution
\[
K_{\varphi\psi}(\ux):=\pD\hD\bigg[\frac{|\ux|^4}{\sigma_m (2-m)(4-m)}\bigg]=\dfrac{(2-m)|\ux|^{-m}\ux_{\psi}\ux_{\varphi}+|\ux|^{2-m}\sum_{i=1}^{m}\psi_{i}\varphi_{i}}{2\sigma_{m}(2-m)},
\]
which is obviously $(\varphi,\psi)$-harmonic in $\R^m\setminus\{0\}$.

The Laplace operator $\Delta_{m}$ is the quintessential example of a (strongly) elliptic operator, while the operator $\hD\pD$ being elliptic, is not strongly elliptic. This can be shown via the following counterexample, which violates the maximum principle.

Let $m$ be even and consider an arbitrary structural set $\varphi=\{\varphi_1,\varphi_2,\dots,\varphi_{m-1},\varphi_m\}$ together with $\psi=\{\varphi_2,\varphi_1,\dots,\varphi_{m},\varphi_{m-1}\}$. Introduce the function 
\[u(\ux)=1-\sum_{i=1}^m x_i^2,\] 
which obviously does vanishes on the boundary of the unit ball $B(0,1)$ of $\R^m$. 

On the other hand we get
\begin{eqnarray*}
\hD\pD u&=&\hD(\varphi_{2}\dfrac{\pa u}{\pa x_{1}}+\varphi_{1}\dfrac{\pa u}{\pa x_{2}}+\cdots +\varphi_{m}\dfrac{\pa u}{\pa x_{m-1}}+\varphi_{m-1}\dfrac{\pa u}{\pa x_{m}})\\&=&\varphi_{1}\varphi_{2}\dfrac{\pa^{2} u}{\pa x_{1}^{2}}+\varphi_{2}\varphi_{1}\dfrac{\pa^{2} u}{\pa x_{2}^{2}}+\cdots+\varphi_{m-1}\varphi_{m}\dfrac{\pa^{2} u}{\pa x_{m-1}^{2}}+\varphi_{m}\varphi_{m-1}\dfrac{\pa^{2} u}{\pa x_{m}^{2}}\\&=&-2\varphi_{1}\varphi_{2}-2\varphi_{2}\varphi_{1}-\cdots-2\varphi_{m-1}\varphi_{m}-2\varphi_{m}\varphi_{m-1}=0.
\end{eqnarray*}
This give $u\in\cH_{\varphi,\psi}(B(0,1),\R_{0,m})$ and non identically zero.    

Of course, the above example mathematically leads to ill-posed formulation of the Dirichlet problem for $(\varphi,\psi)$-harmonic functions in the sense of Hadamard \cite{H}, which is a marked difference between $\cH_{\varphi,\psi}(\Omega,\R_{0,m})$ and $\cH(\Omega,\R_{0,m})$. 

\section{Riemann-Hilbert boundary value problem for $(\varphi,\psi)$-harmonic functions}
In this section a Riemann-Hilbert boundary value problem for $(\varphi,\psi)$-harmonic functions will be discussed. To do so, some basic facts are firstly introduced. 
\begin{theorem}\cite{BDGS}
Let $u\in C^2(\Omega\cup \Gamma,\R_{0,m})$. Then for $\ux\in\Omega$ it holds that
\begin{eqnarray}\label{BP}
u(\ux)=\int_{\Gamma}K_\psi(\uy-\ux)n_\psi(\uy)u(\uy)d\uy+\int_\Gamma K_{\varphi\psi}(\uy-\ux)n_\varphi(\uy)\pD u(\uy)d\uy-\nonumber\\
\int_\Omega K_{\varphi\psi}(\uy-\ux)[\hD\pD u(\uy)]d\uy.
\end{eqnarray}
\end{theorem}
\noindent
The corresponding Teodorescu operator here is given by
\[
\cT_{\varphi\psi}v(\ux):=-\int_\Omega K_{\varphi\psi}(\uy-\ux)v(\uy)d\uy,
\]
which satisfies the relations (see \cite[Theorem 5]{BDGS})
\begin{equation}\label{Trel}
\pD\cT_{\varphi\psi}v(\ux)=\cT_\varphi v(\ux),\,\,\hD\pD\cT_{\varphi\psi}v(\ux)=\begin{cases}v(\ux)&\mbox{if $\ux\in\Omega_+$}\\0&\mbox{if $\ux\in\Omega_-$}.\end{cases}
\end{equation}
\begin{remark}
By means of a more detailed analysis, it may be shown that the above Borel-Pompeiu formula will remain valid under weaker requirements, namely that  $u\in C^2(\Omega,\R_{0,m})\cap C^1(\Omega\cup \Gamma,\R_{0,m})$ and 
\[
\int_\Omega|\hD\pD u(\uy)|d\uy<+\infty,
\] 
which ensure the existence of all the integrals in \eqref{BP}. No consideration will be given in this paper to the problem of finding the most general conditions for the validity of \eqref{BP}, but instead the previous ones are completely sufficient for our purposes.
\end{remark}
In particular, for $u\in\cH_{\varphi,\psi}(\Omega,\R_{0,m} )$, one has in $\Omega$
\begin{equation}\label{CF}
u(\ux)=\int_{\Gamma}K_\psi(\uy-\ux)n_\psi(\uy)u(\uy)d\uy+\int_\Gamma K_{\varphi\psi}(\uy-\ux)n_\varphi(\uy)\pD u(\uy)d\uy,
\end{equation}  
the last formula being a sort of Cauchy representation formula for $(\varphi,\psi)$-harmonic functions.

The Lipschitz class $\Li(1+\alpha,\Gamma)$ is well adapted to define on it a Cauchy transform arising from \eqref{CF}. More precisely, given a Lipschitz data $\bu\in\Li(1+\alpha,\Gamma)$, the Cauchy transform is defined by
\begin{equation}\label{CT}
\cC_{\varphi\psi}\bu (\ux)=\int_{\Gamma}K_\psi(\uy-\ux)n_\psi(\uy)\tilde{u}(\uy)d\uy+\int_\Gamma K_{\varphi\psi}(\uy-\ux)n_\varphi(\uy)\pD \tilde{u }(\uy)d\uy,
\end{equation}  
where $\tilde{u}$ denotes the Whitney extension of $\bu$ in Theorem \ref{Wh}.

Although the Whitney extension is not unique, one should remark that the appearance of ambiguity in the above definition disappears since the values of $\tilde{u}$ and $\pD \tilde{u}$ on $\Gamma$ are completely determined by the collection $\bu=\{u^{(\bj)},\,|\bj|\le 1\}$.

Of course, the function $\cC_{\varphi\psi}\bu$ is by definition $(\varphi,\psi)$-harmonic in $\R^m\setminus\Gamma$. Also, it should be noted that the weakly singularity of the kernel $K_{\varphi\psi}(\uy-\ux)$ implies that the second  integral in \eqref{CT} does not experience a jump when
$\ux$ is crossing the boundary $\Gamma$. This together with the classical Plemelj-Sokhotski formulas applied to the first Cauchy type integral in \eqref{CT} lead to
\begin{equation}\label{Jump}
[\cC_{\varphi\psi}\bu]^+(\ux)-[\cC_{\varphi\psi}\bu]^-(\ux)=\tilde{u}(\ux)=u^{(0)}(\ux),\,\,\,\ux\in\Gamma,
\end{equation} 
where 
\[
[\cC_{\varphi\psi}\bu]^\pm(\ux)=\lim_{\Omega_\pm\ni\uz\to\ux}\cC_{\varphi\psi}\bu(\uz).
\]
\subsection{Riemann-Hilbert boundary value problem}
The Riemann-Hilbert boundary value problem (RH-problem for short) for an unknown $u\in \cH_{\varphi,\psi}(\Omega_+\cup\Omega_-,\R_{0,m})$ is defined
for $\bg=\{g^{(\bj)},\,|\bj|\le 1\} \in \Li(1+\alpha,\Gamma)$ by the system 
\begin{equation}\label{RHS}
\begin{cases}
u^{+}(\ux)-u^{-}(\ux)A=\tilde{g}(\ux)&\mbox{if $\ux\in\Gamma$},\\
[\pD u]^{+}(\ux)-[\pD u]^{-}(\ux)B=\pD\tilde{g}(\ux)&\mbox{if $\ux\in\Gamma$,}\\u(\infty)=\pD u(\infty)=0,
\end{cases}
\end{equation}  
where $A,B$ are two invertible $\R_{0,m}$-valued constants.

\subsection{The smooth case}
Let us denote
\[
g_*(\ux)=\int_\Gamma K_{\varphi\psi}(\uy-\ux)n_\varphi(\uy)\pD \tilde{g }(\uy)d\uy(-1+B^{-1}A)+\tilde{g}(\ux),\,\,\ux\in\Gamma.
\]
Then, the function
\begin{equation}\label{SolSmooth}
u(\ux):=\begin{cases}\cC_\psi g_*(\ux)+\int_{\Gamma}K_{\varphi\psi}(\uy-\ux)n_\varphi(\uy)\pD \tilde{g}(\uy)d\uy, &\mbox{$\ux\in\Omega_+$},\\\\\cC_\psi g_*(\ux)A^{-1}\!+\!\int_\Gamma K_{\varphi\psi}(\uy-\ux)n_\varphi(\uy)\pD\tilde{g}(\uy)d\uy B^{-1},&\mbox{$\ux\in\Omega_-$}\end{cases}
\end{equation}
belongs to $\cH_{\varphi,\psi}(\Omega_+\cup\Omega_-,\R_{0,m})$ and satisfies the boundary conditions in \eqref{RHS}.

The uniqueness of homogeneous RH-problem (\ref{RHS}) is reduced to prove that it has only the null-solution. 

In fact, since $u\in \cH_{\varphi,\psi}(\Omega_+\cup\Omega_-,\R_{0,m})$, the function $\pD u(\ux)$ is $\varphi$-hyperholomorphic in $\Omega_+\cup\Omega_-$ and has no jump through $\Gamma$. This, together with the vanishing condition $\pD u(\infty)=0$ yields $\pD u\equiv 0$ in $\R^m$, which is clear from the combination of classical Painleve and Liouville theorems. The proof is concluded after using the same arguments for the $\psi$-hyperholomorphic function $u$.

We refer the reader also to \cite[p.307]{G16} and \cite{ZG}, where similar problems for standard $\R_{0,m}$-valued harmonic functions in the smooth context are studied.  
\subsection{The fractal case}
There is an essential difference between the fractal case and those studied before. In fact, for a fractal boundary $\Gamma$, the function given by \eqref{SolSmooth} is useless and meaningless as it stands. 

Throughout this subsection we follow \cite{HN} in assuming that $\Gamma$ is $d$-summable for some $m-1<d<m$, i.e. the improper integral
\[
\int_0^1 N_\Gamma(\tau) \, \tau^{d-1} \, d\tau
\]
converges, where $ N_\Gamma(\tau)$ stands for the minimal number of balls of radius $\tau$ needed to cover $\Gamma$.

The notion of a $d$-summable subset was introduced by Harrison \& Norton in \cite{HN}, who showed that if $\Gamma$ has box dimension \cite{F} less than $d$, then is $d$-summable.

The following lemma can be found in \cite[Lemma 2]{HN} and reveals the specific importance of the notion of $d$-summability of the boundary $\Gamma$ of a Jordan domain $\Omega$ in connection with the Whitney decomposition $\cW$ of $\Omega$ by squares $\cQ$ of diameter $|\cQ|$.
\begin{lemma}\cite{HN}
\label{dsum}
If $\Omega$ is a Jordan domain of $\R^{m}$ and its boundary $\Gamma$ is $d$-summable, then the expression $\sum_{\cQ\in{\mathcal{W}}}|\cQ|^d$, called the $d$-sum of the Whitney decomposition of $\Omega$, is finite.
\end{lemma}

\begin{lemma}\label{lp}
Let $\bg\in\Li(1+\alpha,\Gamma)$, then $\hD\pD\tilde{g}\in L^p(\Omega,\R_{0,m})$  for $p=\displaystyle\frac{m-d}{1-\alpha}$.
\end{lemma}
\pf From Theorem \ref{Wh} $(iii)$, we have $|\hD\pD\tilde{g}(\ux)| \leqslant c \, \mbox{\em dist}(\ux,\Gamma)^{\alpha-1}$ for $\ux\in\Omega$. With this in mind the proof follows in a quite analogous way to that of \cite[Lemma 4.1]{AAB}.\qed  
\begin{lemma}\label{continuity}
Let $\bg\in\Li(1+\alpha,\Gamma)$ with $\alpha>\displaystyle \frac{d}{m}$. Then $\cT_\psi[\pD\tilde{g}]$ and $\pD\cT_{\varphi\psi}[\hD\pD\tilde{g}]$ are continuous functions in $\R^m$.
\end{lemma}
\pf Since $\pD\tilde{g}\in C^0(\Omega\cup \Gamma,\R_{0,m})$, then $\pD\tilde{g}$ belongs to $L_p(\Omega,\R_{0,m})$, with $p>m$. Similar analysis  to that in the proof of \cite [Proposition 8.1]{GHS} can be applied to conclude that $\cT_\psi[\pD\tilde{g}]\in C^0(\R^m,\R_{0,m})$.

To prove the continuity of $\pD\cT_{\varphi\psi}[\hD\pD\tilde{g}]$ we use the first identity in \eqref{Trel}:
\[
\pD\cT_{\varphi\psi}[\hD\pD\tilde{g}]=\cT_\varphi[\hD\pD\tilde{g}]. 
\] 
A fine point here is to note that under the condition $\alpha>\displaystyle\frac{d}{m}$ we can applied Lemma \ref{lp} to ensure that $\hD\pD\tilde{g}$ also belongs to $L_p(\Omega,\R_{0,m})$ with $p>m$, and so, the reasoning as described above may be repeated.\qed 
 
The following theorem describes a solution of the RH-problem \eqref{RHS} in fractal context.
\begin{theorem}
Let $\bg\in\Li(1+\alpha,\Gamma)$and let $\Gamma$ be $d$-summable  with $\alpha>\displaystyle\frac{d}{m}$. Then the RH-problem \eqref{RHS} has a solution given by
\begin{equation}\label{SolFractal}
u(\ux):=\begin{cases}g_*(\ux)-\cT_\psi[\pD g_*](\ux)+\cT_\psi[\pD\tilde{g}](\ux)-\cT_{\varphi\psi}[\hD\pD\tilde{g}](\ux),&\mbox{$\ux\in\Omega_+$},\\\\-\cT_\psi[\pD g_*](\ux)A^{-1}+\bigg(\cT_\psi[\pD\tilde{g}](\ux)-\cT_{\varphi\psi}[\hD\pD\tilde{g}](\ux)\bigg)B^{-1},&\mbox{$\ux\in\Omega_-$,}\end{cases}
\end{equation}
where
\[
g_*(\ux):=\bigg(\cT_\psi[\pD\tilde{g}](\ux)-\cT_{\varphi\psi}[\hD\pD\tilde{g}](\ux)\bigg)(-1+B^{-1}A)+\tilde{g}(\ux).
\]
\end{theorem}
\pf Lemma \ref{continuity} shows that $g_*\in C^0(\R^m,\R_{0,m})$. On the other hand 
\[
\pD g_*(\ux)=\bigg(\pD\tilde{g}(\ux)-\cT_{\varphi}[\hD\pD\tilde{g}](\ux)\bigg)(-1+B^{-1}A)+\pD\tilde{g}(\ux)
\]
and clearly we will again have a function of $C^0(\R^m,\R_{0,m})$.

Consequently,  $\cT_\psi[\pD g_*]$ together with the remaining terms in \eqref{SolFractal} all belong to $C^0(\R^m,\R_{0,m})$. The first boundary condition in \eqref{RHS} then follows directly after a direct calculation.

To prove the second condition we use again the identities \eqref{rel1} and \eqref{Trel} to obtain

\begin{equation}\label{partialu}
\pD u(\ux):=\begin{cases}\pD\tilde{g}(\ux)-\cT_{\varphi}[\hD\pD\tilde{g}](\ux),&\mbox{$\ux\in\Omega_+$},\\\\-\cT_{\varphi}[\hD\pD\tilde{g}](\ux)B^{-1},&\mbox{$\ux\in\Omega_-$.}\end{cases}
\end{equation} 
Hence, the second boundary condition is directly deduced. Finally, the $(\varphi,\psi)$-harmonicity of $u$ in $\R^m\setminus\Gamma$ follows immediately acting $\hD$ in \eqref{partialu}.\qed 

If $A=B$, then $g_*=\tilde{g}$ and the solution is considerably simplified to
\begin{equation}\label{SolSimply}
u(\ux):=\begin{cases}\tilde{g}(\ux)-\cT_{\varphi\psi}[\hD\pD\tilde{g}](\ux),&\mbox{$\ux\in\Omega_+$},\\\\-\cT_{\varphi\psi}[\hD\pD\tilde{g}](\ux)B^{-1},&\mbox{$\ux\in\Omega_-$.}\end{cases}
\end{equation}
\begin{remark}
The method followed to prove the uniqueness of solution of the RH-problem breaks down when we drop the smoothness assumption over the boundary of the domain. Some partial evidence support the conjecture that the uniqueness may be ensured under similar conditions to those discussed in \cite[Theorem 4.2]{AAB}. 
\end{remark}

\end{document}